\documentclass{amsart}
\usepackage{amsmath,amssymb,amsthm}
\usepackage[colorlinks,allcolors=black]{hyperref}

\newtheorem{theorem}{Theorem}[section]
\newtheorem{lemma}[theorem]{Lemma}
\newtheorem{proposition}[theorem]{Proposition}
\newtheorem{corollary}[theorem]{Corollary}
\theoremstyle{definition}
\newtheorem{definition}[theorem]{Definition}
\newtheorem{notation}[theorem]{Notation}
\newtheorem{example}[theorem]{Example}
\theoremstyle{remark}
\newtheorem{remark}[theorem]{Remark}

\newcommand{\Max}{\operatorname{Max}}
\newcommand{\Aut}{\operatorname{Aut}}
\newcommand{\Sym}{\operatorname{Sym}}
\newcommand{\lk}{\operatorname{lk}}
\newcommand{\CG}{\mathrm{CG}_J}
\newcommand{\cC}{\mathcal{C}}
\newcommand{\cF}{\mathcal{F}}
\newcommand{\cP}{\mathcal{P}}
\newcommand{\cH}{\mathcal{H}}
\newcommand{\nn}{[n]}
\newcommand{\doi}[1]{\href{https://doi.org/#1}{doi:#1}}
\begin{document}

\title[Minimal comaximal hypergraphs]{Isomorphisms and automorphisms of minimal comaximal hypergraphs of commutative rings}

\author{Elham Mehdi-Nezhad}
\address{Department of Mathematics and Applied Mathematics, University of the Western Cape, Private Bag X17, Bellville 7535, Cape Town, South Africa}
\email{emehdinezhad@uwc.ac.za}


\begin{abstract}
Let $R$ be a commutative ring with exactly $n\ge 2$ maximal ideals and let $k\ge 2$. We study the minimal $k$-uniform comaximal hypergraph $\cH_k(R)$, whose vertices are the nonunits outside the Jacobson radical and whose edges are the irredundant $k$-element comaximal families. For $2\le k\le n$ we show that $\cH_k(R)$ is a blow-up of the hypergraph of irredundant $k$-coverings of $\{1,\dots,n\}$, and that its reduced form reconstructs the maximal ideals up to permutation; for $k>n$ there are no edges. We obtain a general description of $\Aut(\cH_k(R))$ and characterize isomorphisms of such hypergraphs. For $3\le k\le n$ its automorphism group properly contains that of the comaximal graph, and for finite rings with $n\ge 3$ and every $2\le k\le n$, $\cH_k(R)$ has the same isomorphism invariants as the comaximal graph, namely $|J(R)|$ and $R/J(R)$. The case $k=2$ recovers and extends the automorphism results for comaximal graphs obtained previously. We also compare the minimal and non-minimal comaximal hypergraphs and identify a condition under which the latter has the same automorphism group and isomorphism invariants as the comaximal graph. The same method determines the vertex set and the automorphism group of the $k$-zero-divisor hypergraph of a product of integral domains, and shows that for finite products of finite fields this hypergraph determines the ring.
\end{abstract}

\maketitle

\section{Introduction}

All rings are commutative with identity. For a ring $R$, we denote by $U(R)$, $J(R)$ and $\Max(R)$ the group of units, the Jacobson radical and the set of maximal ideals of $R$, respectively. Sharma and Bhatwadekar \cite{SB95} introduced the graph on $R$ in which distinct $a,b$ are adjacent if $Ra+Rb=R$; Maimani, Salimi, Sattari and Yassemi \cite{MSSY08} studied it under the name comaximal graph and concentrated on the subgraph induced on the nonunits. Following \cite{KMR25} we write $\CG(R)$ for the subgraph induced on
\[
V(R)=R\setminus\bigl(U(R)\cup J(R)\bigr).
\]
The automorphism group of $\CG(R)$ was studied in \cite{KMR25}: it is determined there for rings with two maximal ideals \cite[Proposition~3.6]{KMR25} and, for finite rings with three maximal ideals, under the hypothesis that the six strata $T_i$, $I_{ij}$ have pairwise distinct degrees \cite[Corollary~3.10]{KMR25}; for $n\ge3$ maximal ideals, \cite[Theorem~3.9]{KMR25} treats a subgraph under a similar distinctness hypothesis.

Hypergraph analogues of graphs attached to rings usually impose a minimality condition on edges. In the $k$-zero-divisor hypergraph of Eslahchi and Rahimi \cite{ER07} a $k$-set is an edge when its product is zero and no proper product vanishes; see also \cite{HB12} for a modified version. For a finite direct product of $k\ge3$ integral domains, the vertices of the $k$-zero-divisor hypergraph are exactly the elements with one zero component, and this hypergraph is connected with diameter at most $3$ and domination number at most $2$ \cite[Lemma~2.1, Theorem~2.1]{MNR23}; for Artinian rings that are products of $k\ge3$ local rings the bounds are $4$ and $2k$ \cite[Theorem~2.2]{MNR23}. Dominating sets of comaximal graphs were studied in \cite{MNR15}. The comaximal analogue of the $k$-zero-divisor hypergraph is the following.

\begin{definition}\label{def:H}
Let $k\ge 2$. The \emph{minimal comaximal hypergraph} $\cH_k(R)$ is the $k$-uniform hypergraph with vertex set $V(R)$ in which a $k$-subset $\{a_1,\dots,a_k\}\subseteq V(R)$ is an edge if and only if
\[
Ra_1+\dots+Ra_k=R \quad\text{and}\quad \sum_{i\in I}Ra_i\ne R\ \text{ for every proper subset } I\subsetneq\{1,\dots,k\}.
\]
We write $\cH_k^{\circ}(R)$ for the hypergraph obtained by deleting the isolated vertices.
\end{definition}

Since a singleton $\{a\}$ with $a\notin U(R)$ never generates $R$, $\cH_2(R)=\CG(R)$.
Selvakumar and Amritha \cite{SA20,SA19} introduced the \emph{$k$-maximal hypergraph} of a commutative ring, whose edges are $k$-sets $\{a_1,\dots,a_k\}$ with $\sum Ra_i=R$ subject to a minimality condition, and studied its connectedness, diameter, girth, regularity, colourings, planarity and genus. Its edges are minimal comaximal $k$-families as in Definition~\ref{def:H}, so it differs from $\cH_k(R)$ at most in the treatment of vertices that lie in no edge. The isomorphism and automorphism problems developed here are not among the principal results of \cite{SA20,SA19}. A non-minimal comaximal hypergraph, whose hyperedges are the finite subsets $\{x_0,\dots,x_r\}\subseteq R$ with $Rx_0+\dots+Rx_r=R$, was considered by Milo\v{s}evi\'c \cite{Mil15} in the study of independence complexes of comaximal graphs. We study the minimal $k$-uniform version; the $k$-uniform part of the non-minimal hypergraph on $V(R)$ is compared with $\CG(R)$ in Section~\ref{sec:full}.

Throughout, $R$ has exactly $n\ge 2$ maximal ideals $M_1,\dots,M_n$, and $\nn=\{1,\dots,n\}$. For $x\in R$ put
\[
S(x)=\{i: x\in M_i\},\qquad C(x)=\nn\setminus S(x)=\{i: x\notin M_i\}.
\]
We call $C(x)$ the \emph{co-support} of $x$. For $\emptyset\ne C\subsetneq\nn$ let $X(C)=\{x\in V(R): C(x)=C\}$.

The main results are the following.

\begin{enumerate}
\item[(A)] (Structure, Section~\ref{sec:structure}.) $\{a_1,\dots,a_k\}$ is an edge of $\cH_k(R)$ if and only if $C(a_1),\dots,C(a_k)$ form an irredundant covering of $\nn$. Thus $\cH_k(R)$ is a blow-up of the hypergraph $\cC_k(n)$ of irredundant $k$-coverings of $\nn$. A vertex $x$ is isolated if and only if $x$ lies in at most $k-2$ maximal ideals; for $k>n$ there are no edges.
\item[(B)] (Reconstruction and rigidity, Sections~\ref{sec:links}--\ref{sec:rigid}.) Two vertices of $\cH_k(R)$ have the same link if and only if they have the same co-support or are both isolated; so the fibres, and with them the abstract covering hypergraph $\cC_k^{\circ}(n)$, are recovered from $\cH_k(R)$. Every isomorphism $\cC_k^{\circ}(n)\to\cC_k^{\circ}(n')$ is induced by a unique bijection $\nn\to[n']$.
\item[(C)] (Automorphisms, Theorem~\ref{thm:aut}.) For $2\le k\le n$, $\Aut(\cH_k(R))=\Sym(I_k(R))\times\bigl(\bigl(\prod_{|C|\le n-k+1}\Sym(X(C))\bigr)\rtimes\Gamma_k(R)\bigr)$. For finite $R$, $\Gamma_k(R)$ is the group of permutations $\sigma$ of $\nn$ with $|R/M_{\sigma(i)}|=|R/M_i|$ for all $i$; it does not depend on $k$.
\item[(D)] (Isomorphisms, Theorems~\ref{thm:iso} and \ref{thm:finiteiso}.) For finite rings with $n\ge 3$ and $2\le k\le n$: $\cH_k(R)\cong\cH_k(R')$ $\iff$ $\CG(R)\cong\CG(R')$ $\iff$ $|J(R)|=|J(R')|$ and $R/J(R)\cong R'/J(R')$.
\item[(E)] (Comparison, Sections~\ref{sec:compare}--\ref{sec:full}.) $\Aut(\CG(R))$ is a proper subgroup of $\Aut(\cH_k(R))$ for $3\le k\le n$. For finite rings with $n\ge3$ and $2\le k\le n-1$, the hypergraph $\cH_k^{\circ}(R)$ determines $R$ up to the invariants in (D), but $\cH_n^{\circ}(R)$ does not. The non-minimal hypergraph, in which every comaximal $k$-set is an edge, has the same automorphism group and the same isomorphism invariants as $\CG(R)$ whenever $|M_i\setminus J(R)|\ge k$ for all $i$.
\end{enumerate}

\subsection*{Relationship to previous work}
The case $k=2$ is the comaximal graph studied in \cite{KMR25}, where automorphism groups were computed by degree and stratum arguments under distinctness hypotheses. For $k=2$, Theorem~\ref{thm:aut} describes $\Aut(\CG(R))$ for every semilocal ring with at least two maximal ideals, and Theorems~\ref{thm:iso} and \ref{thm:finiteiso} give isomorphism criteria. No distinctness hypothesis is needed: the symmetries between fibres are exactly those coming from $\Gamma_k(R)$ in (C). For $k\ge3$ the isolated vertices give the extra factor $\Sym(I_k(R))$, and Theorem~\ref{thm:compare} computes the index of $\Aut(\CG(R))$ in $\Aut(\cH_k(R))$. In Section~\ref{sec:zd} the same method gives the vertex set (for $k=n$ see \cite[Lemma~2.1]{MNR23}), the automorphism group of the $k$-zero-divisor hypergraph of a product of integral domains, and its isomorphism type for products of finite fields; the domination bound of \cite[Theorem~2.1]{MNR23} is the case $k=n$ of the bound $n-k+2$ in Proposition~\ref{prop:dom}. The automorphism group of a blow-up of a rigid hypergraph has a standard description; the work here lies in identifying the covering hypergraph of the ring and proving its rigidity (Proposition~\ref{prop:rigid}).

Except in Section~\ref{sec:finite} and where stated, $R$ need not be finite: we only assume that $R$ is \emph{semilocal}, i.e.\ $R$ has only finitely many maximal ideals.

A $k$-uniform hypergraph is a pair $H=(V,E)$ with $E$ a set of $k$-subsets of $V$. An isomorphism $H\to H'$ is a bijection $V\to V'$ that maps $E$ onto $E'$. For a set $X$, $\Sym(X)$ is its full symmetric group (also when $X$ is infinite), and $\Sym(n)=\Sym(\nn)$.

\section{Co-supports and the covering hypergraph}\label{sec:structure}

Let $\pi\colon R\to\prod_{i=1}^n R/M_i$ be the canonical map. Since the $M_i$ are pairwise comaximal, $\pi$ is surjective with kernel $\bigcap_i M_i=J(R)$ by the Chinese Remainder Theorem \cite[Proposition~1.10]{AM69}. Clearly $x\in M_i$ if and only if the $i$-th coordinate of $\pi(x)$ is zero. Let $\cP_n$ be the set of nonempty proper subsets of $\nn$.

\begin{lemma}\label{lem:basic}
\begin{enumerate}
\item[(a)] $x\in V(R)$ if and only if $C(x)\in\cP_n$.
\item[(b)] For $a_1,\dots,a_r\in R$, $Ra_1+\dots+Ra_r=R$ if and only if $C(a_1)\cup\dots\cup C(a_r)=\nn$.
\item[(c)] For every $C\in\cP_n$ the fibre $X(C)$ is nonempty and is a union of cosets of $J(R)$, and
\[
|X(C)|=|J(R)|\cdot\prod_{i\in C}\bigl|(R/M_i)^{\times}\bigr| \qquad(\text{cardinal arithmetic}).
\]
\end{enumerate}
\end{lemma}

\begin{proof}
(a) $x\notin U(R)$ iff $x$ lies in some maximal ideal iff $C(x)\ne\nn$; $x\notin J(R)$ iff $x$ avoids some $M_i$ iff $C(x)\ne\emptyset$.

(b) The ideal $\sum_i Ra_i$ is proper iff it lies in some $M_l$, iff $a_1,\dots,a_r\in M_l$ for some $l$, iff some $l$ lies in no $C(a_i)$.

(c) $X(C)=\pi^{-1}(Y_C)$ with $Y_C=\{y: y_i\ne 0\iff i\in C\}$. The set $Y_C$ is nonempty of cardinality $\prod_{i\in C}|(R/M_i)^\times|$, and the fibres of the surjection $\pi$ are the cosets of $J(R)$.
\end{proof}

\begin{definition}
A family $(C_1,\dots,C_k)$ of subsets of $\nn$ is an \emph{irredundant covering} of $\nn$ if $\bigcup_i C_i=\nn$ and no proper subfamily covers $\nn$. Equivalently, it covers $\nn$ and every $C_i$ has a \emph{private element}, i.e.\ an element lying in no $C_j$ with $j\ne i$.
\end{definition}

If $(C_1,\dots,C_k)$ is an irredundant covering and $k\ge 2$, the $C_i$ are pairwise distinct (if $C_1=C_2$, then $C_2$ has no private element) and nonempty (each has a private element).

\begin{lemma}\label{lem:edge}
Let $k\ge 2$ and let $a_1,\dots,a_k\in V(R)$ be distinct. Then $\{a_1,\dots,a_k\}$ is an edge of $\cH_k(R)$ if and only if $(C(a_1),\dots,C(a_k))$ is an irredundant covering of $\nn$. In this case the co-supports $C(a_i)$ are pairwise distinct.
\end{lemma}

\begin{proof}
By Lemma~\ref{lem:basic}(b), applied to the whole family and to every proper subfamily, the two conditions in Definition~\ref{def:H} say exactly that $\bigcup_i C(a_i)=\nn$ and that no proper subfamily covers $\nn$. The last assertion was noted above.
\end{proof}

\begin{definition}
For $k\ge 2$ let $\cC_k(n)$ be the $k$-uniform hypergraph with vertex set $\cP_n$ whose edges are the $k$-subsets $\{C_1,\dots,C_k\}\subseteq\cP_n$ that form irredundant coverings of $\nn$.
\end{definition}

\begin{corollary}\label{cor:blowup}
A $k$-subset $e\subseteq V(R)$ is an edge of $\cH_k(R)$ if and only if $C$ is injective on $e$ and $C(e)$ is an edge of $\cC_k(n)$. Hence $\cH_k(R)$ is obtained from $\cC_k(n)$ by replacing each vertex $C$ by the independent set $X(C)$ and each edge $\{C_1,\dots,C_k\}$ by all transversals of $X(C_1),\dots,X(C_k)$. In particular, any permutation of $V(R)$ that maps each fibre $X(C)$ onto itself is an automorphism of $\cH_k(R)$.
\end{corollary}

\begin{proof}
Immediate from Lemma~\ref{lem:edge}, since edges of $\cH_k(R)$ have pairwise distinct co-supports and the edge condition depends only on co-supports.
\end{proof}

For the rest of the paper put
\[
m=m(n,k)=n-k+1,\qquad \cP_n^{\le m}=\{C\in\cP_n: |C|\le m\}.
\]

\begin{lemma}\label{lem:comb}
Let $k\ge 2$; in (b) and (c) assume moreover $k\le n$.
\begin{enumerate}
\item[(a)] Every member of an irredundant $k$-covering of $\nn$ has at most $m$ elements. If $k>n$, $\nn$ has no irredundant $k$-covering.
\item[(b)] Let $C\in\cP_n$ with $|C|\le m$ and let $B_1,\dots,B_{k-1}$ be a partition of $\nn\setminus C$ into nonempty blocks (such a partition exists). Then $\{C,B_1,\dots,B_{k-1}\}$ is an edge of $\cC_k(n)$.
\item[(c)] Let $k\ge 3$ and let $C\ne D$ in $\cP_n$. Then $C$ and $D$ lie in a common edge of $\cC_k(n)$ if and only if $C$ and $D$ are incomparable under inclusion and $|C\cup D|\le m+1$.
\end{enumerate}
\end{lemma}

\begin{proof}
(a) Let $(C_1,\dots,C_k)$ be an irredundant covering and choose a private element $p_j$ of each $C_j$. The $p_j$ are distinct, so $k\le n$. Moreover $C_i$ contains none of the $k-1$ elements $p_j$, $j\ne i$, so $|C_i|\le n-k+1$.

(b) $|\nn\setminus C|\ge n-m=k-1$, so the partition exists. The family covers $\nn$; every element of $C$ is private to $C$, and every element of $B_j$ is private to $B_j$. The members are pairwise distinct, nonempty and proper (the $B_j$ are disjoint from $C\neq\emptyset$).

(c) Suppose $\{C,D,C_3,\dots,C_k\}$ is an edge. Then $C$ has a private element, so $C\not\subseteq D$; likewise $D\not\subseteq C$. The private elements of $C_3,\dots,C_k$ are $k-2$ distinct elements outside $C\cup D$, so $|C\cup D|\le n-k+2=m+1$.
Conversely, let $C,D$ be incomparable with $|C\cup D|\le m+1$, and let $W=\nn\setminus(C\cup D)$. Then $|W|\ge n-m-1=k-2\ge 1$; partition $W$ into nonempty blocks $B_1,\dots,B_{k-2}$. Pick $c\in C\setminus D$ and $d\in D\setminus C$. In $(C,D,B_1,\dots,B_{k-2})$, $c$ is private to $C$, $d$ is private to $D$ and every element of $B_j$ is private to $B_j$; the family covers $\nn$, and its members are pairwise distinct proper nonempty subsets. Hence it is an edge of $\cC_k(n)$.
\end{proof}

By Lemma~\ref{lem:comb}(a),(b), the non-isolated vertices of $\cC_k(n)$ are exactly the sets in $\cP_n^{\le m}$. Write $\cC_k^{\circ}(n)$ for the hypergraph on $\cP_n^{\le m}$ with the same edges.

\begin{corollary}\label{cor:isolated}
Let $2\le k\le n$. A vertex $x$ of $\cH_k(R)$ is isolated if and only if $|C(x)|\ge n-k+2$, i.e.\ if and only if $x$ belongs to at most $k-2$ maximal ideals. Thus the set $I_k(R)$ of isolated vertices is $\bigcup_{|C|>m}X(C)$; it is empty for $k=2$ and nonempty for $k\ge3$. If $k>n$, $\cH_k(R)$ has no edges.
\end{corollary}

\begin{proof}
By Corollary~\ref{cor:blowup} and nonemptiness of the fibres, $x$ is non-isolated iff $C(x)$ is a non-isolated vertex of $\cC_k(n)$, i.e.\ iff $|C(x)|\le m$. For $k\ge 3$ the fibre $X(\nn\setminus\{1\})$ is nonempty and $n-1\ge m+1$. The last statement is Lemma~\ref{lem:comb}(a).
\end{proof}

\begin{example}\label{ex:kn}
(i) $k=n$. Then $m=1$ and the only edge of $\cC_n(n)$ is $\{\{1\},\dots,\{n\}\}$. Hence $\cH_n^{\circ}(R)$ is the complete $n$-partite $n$-uniform hypergraph with parts
\[
X(\{i\})=\Bigl(\bigcap_{l\ne i}M_l\Bigr)\setminus M_i\qquad(1\le i\le n),
\]
and every element lying in at most $n-2$ maximal ideals is isolated.

(ii) $k=n-1\ge 2$. Then $m=2$. Let $(C_1,\dots,C_{n-1})$ be an irredundant covering with private elements $e_1,\dots,e_{n-1}$ and let $z$ be the remaining element of $\nn$. Since $C_i$ contains no $e_j$ with $j\ne i$, we have $C_i\in\{\{e_i\},\{e_i,z\}\}$, and $z$ lies in at least one $C_i$. Conversely every such family is an irredundant covering. A covering in which $z$ lies in exactly one member $C_i=\{e_i,z\}$ arises from two such descriptions ($e_i$ and $z$ may be interchanged); otherwise the description is unique. Counting pairs $(z,A)$, where $A=\{i: z\in C_i\}\neq\emptyset$, gives
\[
|E(\cC_{n-1}(n))|=n\bigl(2^{n-1}-1\bigr)-\binom n2 .
\]
\end{example}

\section{Links and fibres}\label{sec:links}

For a vertex $x$ of a $k$-uniform hypergraph $H$, the \emph{link} of $x$ is $\lk(x)=\{e\setminus\{x\}: x\in e\in E(H)\}$. Write $x\equiv y$ if $\lk(x)=\lk(y)$. This is an equivalence relation, and every isomorphism $\varphi\colon H\to H'$ satisfies $\lk(\varphi(x))=\varphi(\lk(x))$, so it maps $\equiv$-classes onto $\equiv$-classes. Note that $x\equiv y$ with $x\neq y$ forces $x,y$ to lie in no common edge: if $x,y\in e$, then $e\setminus\{x\}\in\lk(x)$ contains $y$, whereas no member of $\lk(y)$ contains $y$.

\begin{lemma}\label{lem:separate}
Let $2\le k\le n$ and let $C\ne D$ in $\cP_n^{\le m}$ with $C\not\subseteq D$. Then either $C$ and $D$ lie in a common edge of $\cC_k(n)$, or there is a $(k-1)$-subset $\cF\subseteq\cP_n\setminus\{C,D\}$ such that $\cF\cup\{C\}$ is an edge and $\cF\cup\{D\}$ is not.
\end{lemma}

\begin{proof}
Fix $a\in C\setminus D$ and let $B=\nn\setminus C$, so $b:=|B|\ge k-1$. Choose distinct $p_1,\dots,p_{k-2}\in B$ and put
\[
\cF=\bigl\{\{p_1\},\dots,\{p_{k-2}\},\,B\setminus\{p_1,\dots,p_{k-2}\}\bigr\},
\]
a partition of $B$ into $k-1$ nonempty blocks; its blocks have sizes $1$ and $b-k+2$. By Lemma~\ref{lem:comb}(b), $\cF\cup\{C\}$ is an edge. Since $a\notin D$ and $a\notin\bigcup_{G\in\cF}G=B$, the family $\cF\cup\{D\}$ does not cover $\nn$. The blocks are disjoint from $C\ne\emptyset$, so $C\notin\cF$. It remains to treat the case $D\in\cF$, which forces $D\subseteq B$, i.e.\ $C\cap D=\emptyset$.

If $k=2$, then $\cF=\{B\}$ and $D=B$, so $\{C,D\}$ is a partition of $\nn$ into two nonempty sets, which is an edge. If $k\ge3$ and $|C\cup D|\le m+1$, then $C,D$ (disjoint and nonempty, hence incomparable) lie in a common edge by Lemma~\ref{lem:comb}(c). Finally let $k\ge3$ and $|C|+|D|=|C\cup D|\ge m+2$. Then $|D|\ge m+2-|C|=b-k+3$, so $D$ is not the block of size $b-k+2$; and $D$ is not a singleton block, since $b-k+3\ge 2$. Thus $D\notin\cF$.
\end{proof}

\begin{proposition}\label{prop:equiv}
Let $2\le k\le n$ and $x\ne y$ in $V(R)$. Then $x\equiv y$ in $\cH_k(R)$ if and only if either $C(x)=C(y)$, or $x$ and $y$ are both isolated.
\end{proposition}

\begin{proof}
Suppose $C(x)=C(y)$ and let $F\in\lk(x)$. By Lemma~\ref{lem:edge}, $y\notin F$, and replacing $x$ by $y$ does not change the family of co-supports, so $F\in\lk(y)$. By symmetry $\lk(x)=\lk(y)$. If $x,y$ are isolated, both links are empty.

Conversely, let $C(x)=C\neq D=C(y)$. If exactly one of $x,y$ is isolated, the links differ. Suppose neither is isolated; then $C,D\in\cP_n^{\le m}$ by Corollary~\ref{cor:isolated}. As $\equiv$ is symmetric we may assume $C\not\subseteq D$. Apply Lemma~\ref{lem:separate}.
If $\{C,D,C_3,\dots,C_k\}$ is an edge of $\cC_k(n)$, choose $z_i\in X(C_i)$ for $3\le i\le k$; then $\{x,y,z_3,\dots,z_k\}$ is an edge of $\cH_k(R)$ by Corollary~\ref{cor:blowup}, and $\lk(x)\ne\lk(y)$ by the remark preceding Lemma~\ref{lem:separate}.
Otherwise let $\cF$ be as in Lemma~\ref{lem:separate} and choose $z_G\in X(G)$ for each $G\in\cF$. Since $C,D\notin\cF$, the set $F=\{z_G: G\in\cF\}$ does not contain $x$ or $y$. By Corollary~\ref{cor:blowup}, $F\cup\{x\}$ is an edge and $F\cup\{y\}$ is not. Hence $F\in\lk(x)\setminus\lk(y)$.
\end{proof}

\begin{corollary}\label{cor:classes}
For $2\le k\le n$ the $\equiv$-classes of $\cH_k(R)$ are the fibres $X(C)$, $C\in\cP_n^{\le m}$, together with $I_k(R)$ when $k\ge 3$. The number of $\equiv$-classes of non-isolated vertices is
\[
N(n,k)=\sum_{t=1}^{n-k+1}\binom nt ,
\]
which is strictly increasing in $n\ge k$.
\end{corollary}

\begin{proof}
The first assertion is Proposition~\ref{prop:equiv} together with Lemma~\ref{lem:basic}(c). For monotonicity, $N(n+1,k)$ has one more summand than $N(n,k)$, and $\binom{n+1}t\ge\binom nt$.
\end{proof}

Consequently $\cC_k^\circ(n)$ is determined by $\cH_k(R)$ up to isomorphism: its vertices correspond to the non-isolated $\equiv$-classes, and $k$ such classes form an edge if and only if some edge of $\cH_k(R)$ meets each of them.

\section{Rigidity of the covering hypergraph}\label{sec:rigid}

Two vertices of $\cC_k^\circ(n)$ are \emph{co-edge} if they lie in a common edge. Call distinct $C,D\in\cP_n^{\le m}$ \emph{exchangeable} if there is a $(k-1)$-subset $\cF\subseteq\cP_n^{\le m}\setminus\{C,D\}$ such that both $\cF\cup\{C\}$ and $\cF\cup\{D\}$ are edges. Both notions are preserved by isomorphisms.

\begin{lemma}\label{lem:singletons}
Let $3\le k\le n$ and $C\in\cP_n^{\le m}$.
\begin{enumerate}
\item[(a)] If $C$ and $D$ are exchangeable, then $C\cap D\ne\emptyset$.
\item[(b)] If $a\in\nn$ and $D\in\cP_n^{\le m}\setminus\{\{a\}\}$, then $\{a\}$ and $D$ are co-edge if and only if $a\notin D$.
\item[(c)] $C$ is a singleton if and only if every $D\in\cP_n^{\le m}\setminus\{C\}$ that is not co-edge with $C$ is exchangeable with $C$.
\end{enumerate}
\end{lemma}

\begin{proof}
(a) Let $\cF\cup\{C\}$ and $\cF\cup\{D\}$ be edges and $U=\bigcup_{G\in\cF}G$. Since $C$ has a private element in $\cF\cup\{C\}$, $\nn\setminus U\ne\emptyset$; since $\cF\cup\{C\}$ covers $\nn$, $\nn\setminus U\subseteq C$. Likewise $\nn\setminus U\subseteq D$.

(b) If $a\in D$, then $\{a\}\subsetneq D$ and Lemma~\ref{lem:comb}(c) shows they are not co-edge. If $a\notin D$, then $\{a\}$ and $D$ are incomparable and $|D\cup\{a\}|=|D|+1\le m+1$, so they are co-edge by Lemma~\ref{lem:comb}(c).

(c) Let $C=\{a\}$ and let $D$ be not co-edge with $C$. By (b), $a\in D$, and $|D|\ge 2$. Since $|\nn\setminus D|\ge n-m=k-1$, we may choose distinct $p_1,\dots,p_{k-1}\notin D$. Put
\[
G=\nn\setminus\{a,p_1,\dots,p_{k-2}\},\qquad \cF=\bigl\{\{p_1\},\dots,\{p_{k-2}\},G\bigr\}.
\]
Then $p_{k-1}\in G$, $a\notin G$ and $|G|=n-k+1=m$, so $\cF\subseteq\cP_n^{\le m}$; the members of $\cF$ are pairwise distinct, and they differ from $\{a\}$ and from $D$ (the singletons because $|D|\ge2$ and $p_i\neq a$; $G$ because $a\notin G$). The family $\cF\cup\{\{a\}\}$ covers $\nn$, with private elements $a$, $p_1,\dots,p_{k-2}$ and $p_{k-1}$ for $\{a\}$, $\{p_1\},\dots,\{p_{k-2}\}$ and $G$ respectively. The family $\cF\cup\{D\}$ covers $\nn$ because $a\in D$; $a$ is private to $D$, each $p_i$ ($i\le k-2$) is private to $\{p_i\}$ since $p_i\notin D\cup G$, and $p_{k-1}$ is private to $G$ since $p_{k-1}\notin D$. So both are edges, and $C,D$ are exchangeable.

Now let $|C|\ge 2$. The integer $s=m+2-|C|$ satisfies $1\le s\le m$ (as $2\le|C|\le m$) and $s\le n-|C|$ (as $m+2\le n$, i.e.\ $k\ge3$). Choose $D\subseteq\nn\setminus C$ with $|D|=s$. Then $D\in\cP_n^{\le m}$, $|C\cup D|=m+2$, so $C,D$ are not co-edge by Lemma~\ref{lem:comb}(c), and $C\cap D=\emptyset$, so they are not exchangeable by (a).
\end{proof}

\begin{proposition}\label{prop:rigid}
Let $2\le k\le\min(n,n')$ and let $\psi\colon\cC_k^\circ(n)\to\cC_k^\circ(n')$ be an isomorphism. Then $n=n'$ and there is a unique $\sigma\in\Sym(n)$ with $\psi(C)=\sigma(C)$ for all $C\in\cP_n^{\le m}$. In particular $\Aut(\cC_k^\circ(n))\cong\Sym(n)$, acting naturally.
\end{proposition}

\begin{proof}
Every $\sigma\in\Sym(n)$ preserves irredundant coverings, so it induces an automorphism; the action is faithful because singletons are vertices. It remains to show that every isomorphism is induced.

\emph{Case $k=2$.} Here $m=n-1$, $\cP_n^{\le m}=\cP_n$, and $C,D$ are adjacent iff $C\cup D=\nn$. The neighbourhood of $C$ is $N(C)=\{D\in\cP_n: D\supseteq\nn\setminus C\}$. We claim $N(C)\subseteq N(D)$ iff $C\subseteq D$. If $C\subseteq D$, every superset of $\nn\setminus C$ is a superset of $\nn\setminus D$. Conversely, $\nn\setminus C\in N(C)\subseteq N(D)$ gives $\nn\setminus C\supseteq\nn\setminus D$. Hence $\psi$ is an isomorphism of posets $(\cP_n,\subseteq)\to(\cP_{n'},\subseteq)$. It maps minimal elements, the singletons, bijectively onto singletons; so $n=n'$ and $\psi(\{a\})=\{\sigma(a)\}$ for a bijection $\sigma$. Since $C$ is the union of the singletons below it, $\psi(C)=\sigma(C)$.

\emph{Case $k\ge3$.} By Lemma~\ref{lem:singletons}(c), $\psi$ maps singletons onto singletons, so $n=n'$ and $\psi(\{a\})=\{\sigma(a)\}$ for some $\sigma\in\Sym(n)$; also $\psi$ maps non-singletons to non-singletons. Let $|C|\ge2$ and $a\in\nn$. By Lemma~\ref{lem:singletons}(b), $a\in C$ iff $\{a\}$ and $C$ are not co-edge, iff $\{\sigma(a)\}$ and $\psi(C)$ are not co-edge, iff $\sigma(a)\in\psi(C)$. Hence $\psi(C)=\sigma(C)$. Uniqueness of $\sigma$ follows from faithfulness.
\end{proof}

\section{Automorphisms and isomorphisms}\label{sec:main}

\begin{notation}
For $2\le k\le n$ let
\[
\Gamma_k(R)=\bigl\{\sigma\in\Sym(n):\ |X(\sigma(C))|=|X(C)|\ \text{for all } C\in\cP_n^{\le m}\bigr\}.
\]
\end{notation}

For each cardinal $\kappa$ fix a set $\Omega_\kappa$ of cardinality $\kappa$, and for each $C\in\cP_n^{\le m}$ fix a bijection $\beta_C\colon X(C)\to\Omega_{|X(C)|}$. For $\sigma\in\Gamma_k(R)$ define $\tilde\sigma\colon V(R)\to V(R)$ by
\[
\tilde\sigma(x)=\beta_{\sigma(C)}^{-1}\bigl(\beta_C(x)\bigr)\quad (x\in X(C),\ |C|\le m),\qquad \tilde\sigma(x)=x\quad(x\in I_k(R)).
\]
Then $\tilde\sigma$ maps $X(C)$ bijectively onto $X(\sigma(C))$, and a direct computation gives $\widetilde{\sigma\tau}=\tilde\sigma\tilde\tau$. Let $\tilde\Gamma_k=\{\tilde\sigma:\sigma\in\Gamma_k(R)\}$.

\begin{theorem}\label{thm:aut}
Let $R$ be a commutative ring with exactly $n\ge2$ maximal ideals and let $2\le k\le n$. Then
\[
\Aut(\cH_k(R))=\Sym(I_k(R))\times\Bigl(\Bigl(\prod_{C\in\cP_n^{\le m}}\Sym(X(C))\Bigr)\rtimes\tilde\Gamma_k\Bigr),\qquad \tilde\Gamma_k\cong\Gamma_k(R),
\]
and $\Aut(\cH_k^\circ(R))=\bigl(\prod_{C\in\cP_n^{\le m}}\Sym(X(C))\bigr)\rtimes\tilde\Gamma_k$. For $k=2$ this reads $\Aut(\CG(R))=\bigl(\prod_{C\in\cP_n}\Sym(X(C))\bigr)\rtimes\tilde\Gamma_2$. Here $\tilde\Gamma_k$ acts on the normal subgroup $\prod_{C}\Sym(X(C))$ by conjugation, $\tilde\sigma\,\Sym(X(C))\,\tilde\sigma^{-1}=\Sym(X(\sigma(C)))$. If $k>n$, $\cH_k(R)$ has no edges and $\Aut(\cH_k(R))=\Sym(V(R))$.
\end{theorem}

\begin{proof}
Let $\varphi\in\Aut(\cH_k(R))$. By Corollary~\ref{cor:classes}, $\varphi$ maps $I_k(R)$ onto itself and permutes the fibres $X(C)$, $C\in\cP_n^{\le m}$; write $\varphi(X(C))=X(\psi(C))$. If $\{C_1,\dots,C_k\}$ is an edge of $\cC_k^\circ(n)$, pick $x_i\in X(C_i)$; then $\{x_1,\dots,x_k\}$ is an edge, hence so is $\{\varphi(x_1),\dots,\varphi(x_k)\}$, and by Corollary~\ref{cor:blowup} $\{\psi(C_1),\dots,\psi(C_k)\}$ is an edge. Applying this to $\varphi^{-1}$ shows $\psi\in\Aut(\cC_k^\circ(n))$. By Proposition~\ref{prop:rigid} there is a unique $\sigma_\varphi\in\Sym(n)$ with $\psi(C)=\sigma_\varphi(C)$, and $|X(\sigma_\varphi(C))|=|X(C)|$ because $\varphi$ is a bijection; so $\sigma_\varphi\in\Gamma_k(R)$. The map $\varphi\mapsto\sigma_\varphi$ is a homomorphism $\Aut(\cH_k(R))\to\Gamma_k(R)$.

For $\sigma\in\Gamma_k(R)$, $\tilde\sigma$ is an automorphism: it fixes $I_k(R)$, and for a $k$-set $e$ of non-isolated vertices $C(\tilde\sigma(e))=\sigma(C(e))$, so Corollary~\ref{cor:blowup} applies since $\sigma$ preserves irredundant coverings. As $\sigma_{\tilde\sigma}=\sigma$, the homomorphism is onto and $\sigma\mapsto\tilde\sigma$ is a splitting. Its kernel consists of the automorphisms preserving $I_k(R)$ and every fibre $X(C)$, $|C|\le m$, setwise; by Corollary~\ref{cor:blowup} every such permutation is an automorphism (the vertices of $I_k(R)$ are isolated), so the kernel is $\Sym(I_k(R))\times\prod_{|C|\le m}\Sym(X(C))$. Hence $\Aut(\cH_k(R))$ is the semidirect product of this kernel by $\tilde\Gamma_k$. Finally $\Sym(I_k(R))$ commutes with $\prod_{|C|\le m}\Sym(X(C))$ and with $\tilde\Gamma_k$, which act trivially on $I_k(R)$, so it splits off as a direct factor. The statement for $\cH_k^\circ(R)$ is the same argument without $I_k(R)$.
\end{proof}

\begin{theorem}\label{thm:iso}
Let $R,R'$ have exactly $n\ge2$ and $n'\ge2$ maximal ideals, respectively, and let $2\le k\le n$.
\begin{enumerate}
\item[(a)] If $k>n'$, then $\cH_k(R)\not\cong\cH_k(R')$.
\item[(b)] If $k\le n'$, then $\cH_k(R)\cong\cH_k(R')$ if and only if $n=n'$, there is $\sigma\in\Sym(n)$ with $|X'(\sigma(C))|=|X(C)|$ for all $C\in\cP_n^{\le m}$, and $|I_k(R)|=|I_k(R')|$. The same holds for $\cH_k^\circ$ with the last condition omitted.
\end{enumerate}
Here $X'(\cdot)$ denotes the fibres of $R'$ with respect to a fixed numbering of $\Max(R')$.
\end{theorem}

\begin{proof}
(a) $\cH_k(R)$ has edges and $\cH_k(R')$ has none (Corollary~\ref{cor:isolated}).

(b) Let $\varphi$ be an isomorphism. It maps $\equiv$-classes onto $\equiv$-classes and isolated vertices onto isolated vertices, so $N(n,k)=N(n',k)$ and $n=n'$ by Corollary~\ref{cor:classes}, and $\varphi(I_k(R))=I_k(R')$. As in the proof of Theorem~\ref{thm:aut}, $\varphi$ induces an isomorphism $\psi\colon\cC_k^\circ(n)\to\cC_k^\circ(n)$ with $\varphi(X(C))=X'(\psi(C))$, and $\psi$ is induced by some $\sigma\in\Sym(n)$ (Proposition~\ref{prop:rigid}). Conversely, given $\sigma$, choose bijections $X(C)\to X'(\sigma(C))$ for $|C|\le m$ and $I_k(R)\to I_k(R')$; by Corollary~\ref{cor:blowup} their union is an isomorphism.
\end{proof}

\section{Finite rings}\label{sec:finite}

In this section $R$ is finite. Put $j=|J(R)|$ and $q_i=|R/M_i|$, and let
\[
\Gamma(R)=\{\sigma\in\Sym(n): q_{\sigma(i)}=q_i\ \text{for all } i\}\cong\prod_q\Sym(n_q),
\]
where $n_q$ is the number of maximal ideals with residue field of order $q$.

\begin{lemma}\label{lem:finite}
$|X(C)|=j\prod_{i\in C}(q_i-1)$ for every $C\in\cP_n$, and $\Gamma_k(R)=\Gamma(R)$ for every $2\le k\le n$. Moreover
\[
|V(R)|=j\Bigl(\prod_{i=1}^n q_i-\prod_{i=1}^n(q_i-1)-1\Bigr).
\]
\end{lemma}

\begin{proof}
The formula for $|X(C)|$ is Lemma~\ref{lem:basic}(c). If $\sigma\in\Gamma(R)$, the formula gives $|X(\sigma(C))|=|X(C)|$ for all $C$. If $\sigma\in\Gamma_k(R)$, apply the defining condition to $C=\{i\}\in\cP_n^{\le m}$: $j(q_{\sigma(i)}-1)=j(q_i-1)$. For the last formula, $|R|=j\prod q_i$, $|U(R)|=j\prod(q_i-1)$ (an element is a unit iff all coordinates of $\pi(x)$ are nonzero), and $U(R)\cap J(R)=\emptyset$.
\end{proof}

\begin{corollary}\label{cor:order}
For $2\le k\le n$,
\[
|\Aut(\cH_k(R))|=|I_k(R)|!\cdot\prod_{C\in\cP_n^{\le m}}\Bigl(j\prod_{i\in C}(q_i-1)\Bigr)!\cdot\prod_q n_q!\,.
\]
In particular every automorphism of $\cH_k(R)$ fixes every non-isolated fibre setwise if and only if the residue fields of $R$ have pairwise distinct orders.
\end{corollary}

\begin{example}\label{ex:z30}
(i) $R=\mathbb Z_{30}$, $M_1=(2)$, $M_2=(3)$, $M_3=(5)$. Then $j=1$, $(q_1,q_2,q_3)=(2,3,5)$ and $\Gamma(R)=1$. For $k=2$, Theorem~\ref{thm:aut} gives $\Aut(\CG(\mathbb Z_{30}))=\prod_{C\in\cP_3}\Sym(X(C))\cong S_8\times S_4\times S_2\times S_4\times S_2$, in agreement with \cite[Example~3.11]{KMR25}; two pairs of fibres have equal sizes, but since $\Gamma(R)=1$ no automorphism interchanges them. For $k=3$, $X(\{1\})=\{15\}$, $X(\{2\})=\{10,20\}$, $X(\{3\})=\{6,12,18,24\}$, and $I_3(\mathbb Z_{30})$ consists of the $14$ elements lying in exactly one maximal ideal. Thus
\[
\Aut(\cH_3(\mathbb Z_{30}))\cong S_{14}\times S_4\times S_2 .
\]
(ii) $R=\mathbb F_2^{\,n}$. All fibres are singletons and $\Gamma(R)=\Sym(n)$, so $\Aut(\cH_k(\mathbb F_2^{\,n}))\cong\Sym\bigl(\sum_{t=m+1}^{n-1}\binom nt\bigr)\times\Sym(n)$; for $k=2$ this is $\Sym(n)$.
\end{example}

\begin{theorem}\label{thm:finiteiso}
Let $R$ and $R'$ be finite rings, let $R$ have $n\ge 3$ maximal ideals, and let $2\le k\le n$. The following are equivalent:
\begin{enumerate}
\item[(i)] $\cH_k(R)\cong\cH_k(R')$;
\item[(ii)] $\CG(R)\cong\CG(R')$;
\item[(iii)] $|J(R)|=|J(R')|$ and $R/J(R)\cong R'/J(R')$;
\item[(iv)] $R'$ has $n$ maximal ideals, $|J(R)|=|J(R')|$, and the multisets $\{|R/M|: M\in\Max(R)\}$ and $\{|R'/M'|: M'\in\Max(R')\}$ coincide.
\end{enumerate}
If moreover $k\le n-1$, they are also equivalent to
\begin{enumerate}
\item[(v)] $\cH_k^\circ(R)\cong\cH_k^\circ(R')$.
\end{enumerate}
\end{theorem}

\begin{proof}
(iii)$\Leftrightarrow$(iv): $R/J(R)\cong\prod_i R/M_i$ and finite fields of equal order are isomorphic; conversely the factors of a finite product of fields are determined up to isomorphism by the ring (they are its quotients by maximal ideals).

(iv)$\Rightarrow$(i),(ii),(v): choose $\sigma$ with $q'_{\sigma(i)}=q_i$. By Lemma~\ref{lem:finite}, $|X'(\sigma(C))|=|X(C)|$ for all $C\in\cP_n$, so a union of bijections $X(C)\to X'(\sigma(C))$ is an isomorphism $\cH_k(R)\to\cH_k(R')$ for every $k$ (Corollary~\ref{cor:blowup}); it restricts to $\cH_k^\circ$.

(i)$\Rightarrow$(iv) and, for $k\le n-1$, (v)$\Rightarrow$(iv). Since $\cH_k(R)$ has edges, $V(R')\ne\emptyset$, so $R'$ is not local; by Theorem~\ref{thm:iso}(a) it has $n'\ge k$ maximal ideals, and by Theorem~\ref{thm:iso}(b) $n'=n$ and there is $\sigma\in\Sym(n)$ with
\begin{equation}\label{eq:sizes}
j\prod_{i\in C}(q_i-1)=j'\prod_{i\in C}(q'_{\sigma(i)}-1)\qquad(C\in\cP_n^{\le m}).
\end{equation}
Suppose first $k\le n-1$, so $m\ge2$. Taking $C=\{i\}$ and $C=\{i,l\}$ with $l\ne i$ (possible as $n\ge3$) and dividing, we get $q_l-1=q'_{\sigma(l)}-1$ for all $l$; then $C=\{l\}$ gives $j=j'$. This uses only non-isolated vertices, so it applies to (v) as well.
Suppose now $k=n$. Then \eqref{eq:sizes} says that the multisets $\{a_i\}$ and $\{a'_i\}$, $a_i=j(q_i-1)$, coincide, and $|I_n(R)|=|I_n(R')|$, hence $|V(R)|=|V(R')|$. Let $e_r$ denote the $r$-th elementary symmetric polynomial in $a_1,\dots,a_n$. Since $q_i=(a_i+j)/j$, Lemma~\ref{lem:finite} gives
\[
|V(R)|=\frac{\prod_i(a_i+j)-\prod_i a_i}{j^{n-1}}-j=\sum_{r=1}^{n-1}\frac{e_r}{j^{\,r-1}} .
\]
For fixed $a_1,\dots,a_n\ge1$ and $n\ge3$ the right-hand side is strictly decreasing in $j>0$, because $e_2>0$. Hence $j=j'$, and then $q'_{\sigma(i)}=q_i$.

(ii)$\Rightarrow$(iv) is the case $k=2$ of (i)$\Rightarrow$(iv).
\end{proof}

\begin{remark}\label{rem:n2}
(a) For $n=2$ only $k=2$ is possible, $\CG(R)$ is the complete bipartite graph $K_{a_1,a_2}$ with $a_i=j(q_i-1)$ \cite[Theorem~2.2]{MSSY08}, and $\CG(R)\cong\CG(R')$ iff $R'$ has two maximal ideals and $\{a_1,a_2\}=\{a'_1,a'_2\}$. This does not determine $j$: $\CG(\mathbb Z_4\times\mathbb Z_2)\cong\CG(\mathbb F_3\times\mathbb F_3)\cong K_{2,2}$.

(b) The hypothesis $k\le n-1$ in (v) cannot be dropped. By Example~\ref{ex:kn}(i), $\cH_3^\circ(R)$ is the complete $3$-partite $3$-uniform hypergraph with parts of sizes $j(q_i-1)$. For $R=\mathbb F_3^{\,3}$ and $R'=\mathbb Z_4\times\mathbb F_2\times\mathbb F_2$ all parts have size $2$, so $\cH_3^\circ(R)\cong\cH_3^\circ(R')$, whereas $|V(R)|=18\ne 12=|V(R')|$; thus $\cH_3(R)\not\cong\cH_3(R')$ and $\CG(R)\not\cong\CG(R')$.
\end{remark}

\begin{remark}
Every ring automorphism $f$ of $R$ maps $V(R)$ onto itself and preserves comaximality of families, so it restricts to an automorphism of $\cH_k(R)$. The restriction map $\Aut(R)\to\Aut(\cH_k(R))$ is injective (we only use $n\ge2$). Indeed, suppose $f$ fixes $V(R)$ pointwise. Let $u\in U(R)$. Since $M_1+M_2=R$, write $u=m_1+m_2$ with $m_i\in M_i$. If $m_1\in J(R)\subseteq M_2$, then $u\in M_2$, which is absurd; so $m_1\notin J(R)$, and likewise $m_2\notin J(R)$. Thus $m_1,m_2\in V(R)$ and $f(u)=u$. If $x\in J(R)$, then $1-x\in U(R)$, so $f(x)=x$. Hence $f=1_R$. For finite rings this is the argument of \cite[Theorem~1.2]{KMR25}.
\end{remark}

\section{Comparison with the comaximal graph}\label{sec:compare}

\begin{theorem}\label{thm:compare}
Let $R$ have exactly $n\ge2$ maximal ideals and $2\le k\le n$. Then $\Aut(\CG(R))\le\Aut(\cH_k(R))$ as subgroups of $\Sym(V(R))$, with equality if and only if $k=2$. If $R$ is finite, the index is the multinomial coefficient
\[
[\Aut(\cH_k(R)):\Aut(\CG(R))]=\frac{|I_k(R)|!}{\prod_{|C|>m}|X(C)|!}.
\]
\end{theorem}

\begin{proof}
Let $\varphi\in\Aut(\CG(R))$. By Theorem~\ref{thm:aut} with $k=2$ there is $\sigma\in\Sym(n)$ with $C(\varphi(x))=\sigma(C(x))$ for all $x\in V(R)$. Since $\sigma$ preserves irredundant coverings, Corollary~\ref{cor:blowup} shows $\varphi\in\Aut(\cH_k(R))$.

Let $k\ge3$; then $n\ge 3$. The sets $C=\nn\setminus\{1\}$ and $D=\nn\setminus\{2\}$ have $n-1\ge m+1$ elements, so $X(C)\cup X(D)\subseteq I_k(R)$. Pick $x\in X(C)$, $y\in X(D)$; the transposition $\tau=(x\,y)$ lies in $\Sym(I_k(R))\le\Aut(\cH_k(R))$. Suppose $\tau\in\Aut(\CG(R))$, with associated $\sigma$. Every fibre $X(E)$ with $E\notin\{C,D\}$ is fixed pointwise by $\tau$, so $\sigma(E)=E$; in particular $\sigma$ fixes all singletons (as $|C|=|D|=n-1\ge2$), so $\sigma=1$ and $\tau$ preserves $X(C)$, contradicting $\tau(x)=y$. The index formula follows from Corollary~\ref{cor:order} and its case $k=2$.
\end{proof}

Thus, for finite rings with $n\ge3$, $\cH_k(R)$ ($2\le k\le n$) and $\cH_k^\circ(R)$ ($2\le k\le n-1$) have the same isomorphism invariants as $\CG(R)$, while $\cH_n^\circ(R)$ has fewer (Remark~\ref{rem:n2}(b)). For $k\ge3$, $\cH_k(R)$ does not determine the co-supports of the vertices in $I_k(R)$; on the remaining vertices it determines them up to $\Gamma(R)$, as $\CG(R)$ does.

\subsection*{Domination}
A set $D$ of vertices of a hypergraph \emph{dominates} if every vertex outside $D$ lies in a common edge with some vertex of $D$ \cite[Definition~1.3]{MNR23}; $\gamma$ denotes the minimum size of a dominating set. For $k=2$ dominating sets of comaximal graphs were studied in \cite{MNR15}.

\begin{proposition}\label{prop:dom}
Let $3\le k\le n$. Then $\gamma(\cH_k^\circ(R))\le n-k+2$. More precisely, if $x_a\in X(\{a\})$ for $1\le a\le n-k+2$, then $\{x_1,\dots,x_{n-k+2}\}$ dominates $\cH_k^\circ(R)$.
\end{proposition}

\begin{proof}
Put $r=m+1=n-k+2$; then $2\le r\le n$. Let $v$ be a non-isolated vertex outside $\{x_1,\dots,x_r\}$ and $D=C(v)$. Since $|D|\le m<r$, some $a\le r$ satisfies $a\notin D$. Then $D\ne\{a\}$, so $\{a\}$ and $D$ lie in a common edge of $\cC_k(n)$ by Lemma~\ref{lem:singletons}(b), so $v$ and $x_a$ lie in a common edge of $\cH_k(R)$ by Corollary~\ref{cor:blowup}.
\end{proof}

For $k=n$ the bound is $2$, which is the comaximal counterpart of \cite[Theorem~2.1]{MNR23}. The bound is not sharp in general: for $R=\mathbb F_2^{\,7}$ and $k=3$, one checks with Lemma~\ref{lem:comb}(c) that the three vertices with co-supports $\{1,2,3\}$, $\{1,4,5\}$ and $\{2,4,6\}$ dominate $\cH_3^\circ(R)$, whereas $n-k+2=6$.

\section{The non-minimal hypergraph}\label{sec:full}

For $k\ge2$ let $\mathcal F_k(R)$ be the $k$-uniform hypergraph on $V(R)$ in which every $k$-set $\{a_1,\dots,a_k\}$ with $\sum Ra_i=R$ is an edge. Thus $\mathcal F_2(R)=\CG(R)$, and here $k$ may exceed $n$. Equivalently, $\mathcal F_k(R)$ is the $k$-uniform part of the restriction to $V(R)$ of the comaximal hypergraph of \cite{Mil15}.

\begin{lemma}\label{lem:Fdet}
A $k$-subset of $V(R)$ is an edge of $\mathcal F_k(R)$ if and only if it contains an edge of $\cH_l(R)$ for some $2\le l\le\min(k,n)$.
\end{lemma}

\begin{proof}
A comaximal family of nonunits contains a minimal comaximal subfamily, of some size $l\ge2$, and $l\le n$ by Lemma~\ref{lem:comb}(a). The converse is clear.
\end{proof}

\begin{theorem}\label{thm:full}
Let $k\ge2$ and suppose $|M_i\setminus J(R)|\ge k$ for every $i$.
\begin{enumerate}
\item[(a)] For distinct $x,y\in V(R)$: $Rx+Ry=R$ if and only if every $k$-subset of $V(R)$ containing $x$ and $y$ is an edge of $\mathcal F_k(R)$.
\item[(b)] $\Aut(\mathcal F_k(R))=\Aut(\CG(R))$.
\item[(c)] If $R'$ also satisfies the hypothesis, then a bijection $V(R)\to V(R')$ is an isomorphism $\mathcal F_k(R)\to\mathcal F_k(R')$ if and only if it is an isomorphism $\CG(R)\to\CG(R')$.
\end{enumerate}
For finite $R$ the hypothesis reads $j\bigl(\prod_{l\ne i}q_l-1\bigr)\ge k$ for all $i$. It holds whenever $n\ge3$ and $k\le 2^{n-1}-1$, in particular for $2\le k\le n$.
\end{theorem}

\begin{proof}
(a) If $Rx+Ry=R$, every family containing $x,y$ is comaximal. If $Rx+Ry\ne R$, then $x,y\in M_i$ for some $i$, and $x,y\in M_i\setminus J(R)$. Choose $k-2$ further elements of $M_i\setminus J(R)$; the resulting $k$-set lies in $M_i$, so it is not an edge.

(c) By (a), applied to $R$ and $R'$, an isomorphism of the hypergraphs preserves adjacency in both directions. Conversely, an isomorphism $\varphi\colon\CG(R)\to\CG(R')$ maps $\equiv$-classes to $\equiv$-classes; arguing as in the proof of Theorem~\ref{thm:iso} with $k=2$, we get $n=n'$ and $\sigma\in\Sym(n)$ with $C(\varphi(x))=\sigma(C(x))$. By Lemma~\ref{lem:basic}(b), $\varphi$ preserves comaximality of arbitrary families.

(b) is (c) with $R'=R$.

For finite $R$, $M_i/J(R)\cong\prod_{l\ne i}R/M_l$, so $|M_i\setminus J(R)|=j(\prod_{l\ne i}q_l-1)\ge 2^{n-1}-1$, and $2^{n-1}-1\ge n$ for $n\ge3$.
\end{proof}

\begin{example}
The size hypothesis is needed. Let $R=\mathbb F_2^{\,3}$ and $k=4$. Then $|V(R)|=6$ and $|M_i\setminus J(R)|=3$, so no $4$-subset of $V(R)$ lies in a maximal ideal, and $\mathcal F_4(R)$ is the complete $4$-uniform hypergraph on six vertices. Hence $\Aut(\mathcal F_4(R))\cong S_6$, whereas $\Aut(\CG(R))\cong S_3$ by Example~\ref{ex:z30}(ii).
\end{example}

\section{The \texorpdfstring{$k$}{k}-zero-divisor hypergraph of a product of domains}\label{sec:zd}

The combinatorial part of the paper (Sections~\ref{sec:links}--\ref{sec:rigid} and the proofs of Theorems~\ref{thm:aut} and \ref{thm:iso}) uses only the fact that the hypergraph is a blow-up of $\cC_k^\circ(n)$ with nonempty fibres. We record the consequence for the $k$-zero-divisor hypergraph of \cite{ER07}.

Let $T=D_1\times\dots\times D_n$ be a product of $n\ge2$ integral domains. Following \cite{ER07} (see \cite[Definitions~1.1, 1.2]{MNR23}), $\cH^z_k(T)$ is the $k$-uniform hypergraph whose vertex set $Z(T,k)$ consists of the elements lying in some edge, where a $k$-set $\{x_1,\dots,x_k\}$ of distinct nonzero nonunits is an edge if $x_1\cdots x_k=0$ and no product of $k-1$ of them is $0$. For $x\in T$ let $Z(x)=\{i: x_i=0\}$ and, for $C\in\cP_n$, $Y(C)=\{x\in T: Z(x)=C\}$; then $Y(C)\ne\emptyset$ and $|Y(C)|=\prod_{i\notin C}|D_i\setminus\{0\}|$.

\begin{proposition}\label{prop:zd}
Let $2\le k\le n$ and $m=n-k+1$.
\begin{enumerate}
\item[(a)] A $k$-set of distinct nonzero nonunits $\{x_1,\dots,x_k\}$ is an edge of $\cH^z_k(T)$ if and only if $(Z(x_1),\dots,Z(x_k))$ is an irredundant covering of $\nn$. Hence $Z(T,k)=\{x\in T: 1\le|Z(x)|\le m\}$, and $\cH^z_k(T)$ is the blow-up of $\cC_k^\circ(n)$ with fibres $Y(C)$, $C\in\cP_n^{\le m}$. For $k=n$ this is \cite[Lemma~2.1]{MNR23}.
\item[(b)] $\Aut(\cH^z_k(T))=\bigl(\prod_{C\in\cP_n^{\le m}}\Sym(Y(C))\bigr)\rtimes\tilde\Gamma^z$, where $\Gamma^z=\{\sigma\in\Sym(n): |Y(\sigma(C))|=|Y(C)|\ \text{for all } C\in\cP_n^{\le m}\}$.
\item[(c)] If $k\ge3$, then $\gamma(\cH^z_k(T))\le n-k+2$.
\item[(d)] Let $T$ and $T'$ be finite products of $n$ and $n'$ finite fields, with $k\le\min(n,n')$. Then $\cH^z_k(T)\cong\cH^z_k(T')$ if and only if $T\cong T'$. Moreover $\Gamma^z$ is the group of permutations preserving the orders of the factors.
\end{enumerate}
\end{proposition}

\begin{proof}
(a) Since each $D_i$ is a domain, $x_1\cdots x_r=0$ if and only if $Z(x_1)\cup\dots\cup Z(x_r)=\nn$; a nonzero nonunit has $Z(x)\neq\nn$, and a nonunit $x$ with $Z(x)=\emptyset$ cannot occur in an edge, since removing it leaves a covering. Thus the edges are exactly the $k$-sets whose zero sets form an irredundant covering, which then consists of pairwise distinct members of $\cP_n$. The vertex set follows from Lemma~\ref{lem:comb}(a),(b), and for $k=n$ it consists of the elements with exactly one zero component.

(b),(c) The proofs of Proposition~\ref{prop:equiv}, Theorem~\ref{thm:aut} and Proposition~\ref{prop:dom} use only the blow-up structure and the nonemptiness of the fibres, and apply with $Z$ in place of $C$ (there are no isolated vertices).

(d) Let $T=\prod_i\mathbb F_{q_i}$, $c_i=q_i-1$ and $P=\prod_i c_i$, so $|Y(C)|=P/\prod_{i\in C}c_i$. An isomorphism forces $n=n'$ (count $\equiv$-classes, Corollary~\ref{cor:classes}) and, by Proposition~\ref{prop:rigid}, a permutation $\sigma$ with $P/c_a=P'/c'_{\sigma(a)}$ for every singleton $\{a\}\in\cP_n^{\le m}$. Multiplying over $a$ gives $P^{n-1}=P'^{\,n-1}$, so $P=P'$ and $c'_{\sigma(a)}=c_a$ for all $a$; hence $T\cong T'$. The converse is clear. The same computation with $T'=T$ identifies $\Gamma^z$.
\end{proof}

By (a), for a finite reduced ring $T$ both $\cH_k^\circ(T)$ and $\cH^z_k(T)$ are blow-ups of the same covering hypergraph $\cC_k^\circ(n)$, with fibre sizes $\prod_{i\in C}(q_i-1)$ and $\prod_{i\notin C}(q_i-1)$ respectively. For $T=\mathbb F_2^{\,n}$ all fibres are singletons and $x\mapsto 1-x$, which interchanges zero sets and co-supports, is an isomorphism $\cH_k^\circ(T)\to\cH^z_k(T)$.

\begin{remark}
If $R$ has infinitely many maximal ideals, co-supports are subsets of an infinite set and the fibres $X(C)$ need not be nonempty for all $C$ (the Chinese Remainder Theorem only gives finite patterns). The arguments above use nonemptiness of all fibres in an essential way.
\end{remark}

\end{document}